\documentclass[12pt,twoside]{amsart} 

\usepackage{amssymb,latexsym,bbm,graphicx,epsfig,epic,eepic,oldgerm,psfrag,lipsum,comment}

\usepackage{amsmath,amsfonts,amscd,mathrsfs,amsthm}
\usepackage{tikz,wasysym}         \usetikzlibrary{trees}
\usepackage{tikz-cd}
\usepackage{marvosym} 
\usepackage{hyperref} 
\usepackage{a4wide}

\usepackage{xypic} 
\input xy
\xyoption{all}

\usepackage[utf8]{inputenc}
\usepackage[english]{babel}            
\usepackage{xcolor}

\definecolor{amber}{rgb}{1.0, 0.75, 0.0}
\definecolor{banan}{rgb}{1.0, 0.88, 0.21}
\definecolor{gold}{rgb}{1.0, 0.84, 0.0}                
\definecolor{goldenpoppy}{rgb}{0.99, 0.76, 0.0}

\theoremstyle{plain}
\newtheorem{theorem}{Theorem}[section]
\newtheorem{lemma}[theorem]{Lemma}                           
\newtheorem{proposition}[theorem]{Proposition}
\newtheorem{corollary}[theorem]{Corollary}
\newtheorem*{remark*}{Remark}
\newtheorem*{remarks*}{Remarks}
\newtheorem{remark}[theorem]{Remark}
\newtheorem{remarks}[theorem]{Remarks}

\newtheorem*{example*}{Example}
\newtheorem*{examples*}{Examples}
\newtheorem{definition}[theorem]{Definition}
\newtheorem*{definition*}{Definition}

\newtheorem*{question*}{Question}

\newtheorem{conjecture}[theorem]{Conjecture}

\numberwithin{figure}{section}
\numberwithin{equation}{section}

\newcommand{\proofend}{\hspace*{\fill} $\square$\\}

\def\Gr{\operatorname{Gr}}

\newenvironment{brsm}{
  \bigl[ \begin{smallmatrix} }{
  \end{smallmatrix} \bigr]}

\def\1{\:\!}
\def\2{\;\!}
\def\m{\medskip}

\def\aff{\operatorname{aff}}

\def\SW{\operatorname{SW}}
\def\Sympc0{\operatorname{Symp^c_0}}

\def\Int{\operatorname{Int}}

\def\ess{\operatorname{ess}}

\def\supp{\operatorname{supp}}

\def\Cliff{\operatorname{Cliff}}

\def\G{\operatorname{G}}
\def\FS{\operatorname{FS}}

\def\area{\operatorname{area}}

\def\ga{\alpha}

\def\gf{\varphi}

\def\T{\operatorname{T}}

\def\cj{{\mathcal J}}

\def\cl{{\mathcal L}}
\def\cm{{\mathcal M}}

\def\CC{\mathbb{C}}
\def\DD{\mathbb{D}}

\def\NN{\mathbb{N}}

\def\QQ{\mathbb{Q}}
\def\RR{\mathbb{R}}

\def\ZZ{\mathbb{Z}}

\def\R{\operatorname{\mathbb{R}}}
\def\RP{\operatorname{\mathbb{R}P}}
\def\CP{\operatorname{\mathbb{C}P}}

\def\pp{\partial}

\def\fm{{\mathfrak m}}
\def\fM{{\mathfrak M}}

\def\b{\bigskip}
\def\m{\medskip}

\def\proof{\noindent {\it Proof. \;}}

\begin{document}

\title[]{Pinwheels as Lagrangian barriers}


\author[]{Jo\'e Brendel}  
\thanks{JB and FS partially supported by SNF grant 200020-144432/1.}
\address{Jo\'e Brendel,
Institut de Math\'ematiques,
Universit\'e de Neuch\^atel}
\email{joe.brendel@unine.ch}

\author[]{Felix Schlenk}  
\address{Felix Schlenk,
Institut de Math\'ematiques,
Universit\'e de Neuch\^atel}
\email{schlenk@unine.ch}

\keywords{symplectic embedding, Lagrangian submanifold, Lagrangian barrier, pinwheel, Markov number, Lagrange spectrum}

\date{\today}

\begin{abstract}
The complex projective plane~$\CP^2$ contains certain Lagrangian CW-complexes called \emph{pinwheels}, 
which have interesting rigidity properties related to solutions of the Markov equation, see for example~\cite{ES18}. 
We compute the Gromov width of the complement of pinwheels and show that it is strictly smaller than the Gromov width of~$\CP^2$, 
meaning that pinwheels are \emph{Lagrangian barriers}\/ in the sense of Biran~\cite{Bi01}. 
The accumulation points of the set of these Gromov widths are in a simple bijection with the 
Lagrange spectrum below~$3$.
\end{abstract}

\maketitle

\tableofcontents

\section{Introduction and main results}

\subsection{The problem}

One measure for the size of a symplectic manifold $(M,\omega)$ 
is the Gromov width
$$
w_{\G} (M,\omega)  \,=\, \sup \left\{ a \mid B^{2n}(a) \mbox{ symplectically embeds into $(M,\omega)$} \right\} ,
$$
where $B^{2n}(a)$ denotes the open ball of capacity $a = \pi r^2$ equipped with the standard symplectic form. 
In some cases, the Gromov width is much smaller than the obstruction given by volume, see~\cite{Gr85}.

Biran~\cite{Bi01} discovered that there are Lagrangian submanifolds such that~$w_{\G}(M\setminus L, \omega) < w_{\G}(M,\omega)$. These are called \emph{Lagrangian barriers}. Our main result is the following. 

\begin{theorem}  \label{thm:grwidth}
The Gromov width of the complement of a Markov pinwheel~$L_p^{\fm} \subset \CP^2$ is equal to
\begin{equation}
\label{eq:grwidth}
w_{\G} \left(\CP^2 \setminus L_p^{\fm} \right)
\,=\,  \frac{bc}{p},
\end{equation}
where~$(p,b,c)$ is the unique Markov triple with~$p \geqslant b \geqslant c$ which can be obtained from~$\fm$ by mutations preserving~$p$.
\end{theorem}

In particular, since~$\frac{bc}{p} < 1$ for all Markov triples~$(p,b,c) \neq (1,1,1)$, 
Markov pinwheels are Lagrangian barriers. 
Before introducing the relevant terminology used in Theorem~\ref{thm:grwidth}, let us briefly motivate this result. 

As a simple example, take the round sphere $S^2$ of area~$1$,
and let $C_a$ be an embedded circle enclosing area $a \leqslant \frac 12$. 
Then $w_{\G}(S^2) = 1$ while $w_{\G} (S^2 \setminus C_a) = 1-a$.
Hence $C_a$ is the stronger a Lagrangian barrier the larger $a$ is.
In the extreme case $a= \frac 12$, the Lagrangian~$C_a$ is \emph{monotone}. 
Recall that a Lagrangian is called monotone if the corresponding area and Maslov homomorphisms are positively proportional. 
From now on, we restrict our attention to~$M=\CP^2$, which contains the Clifford torus~$T_{\Cliff}$ and the real projective plane~$\RP^2$ 
as monotone Lagrangians. 
Both can be viewed as four-dimensional generalizations of the equatorial circle in~$S^2$, 
and their relative Gromov widths are $w_{\G}(\CP^2 \setminus \RP^2) = \frac 12$ (see~\cite{Bi01}) 
and $w_{\G} \bigl(\CP^2 \setminus \T_{\Cliff} \bigr) = \frac 23$ (see~\cite[Corollary 1.2.5]{BiCo09}). There seem to be no known examples of monotone Lagrangian surfaces in~$\CP^2$ for which
the quotient~$w_{\G}(M \setminus L) / w_{\G}(M)$ is not $\frac 12$ or~$\frac 23$, and the only other known value for this quotient in a closed symplectic four-manifold is~$1$
(the Clifford torus in $S^2 \times S^2$ is not a barrier).
Although Markov pinwheels are not smooth, they are monotone in an appropriate sense -- see the appendix. 
Theorem~\ref{thm:grwidth} thus illustrates that many more quotients of Gromov widths can be realized for monotone Lagrangians if one allows them to be slightly singular. In Section~\ref{s:numbers} we study the set of numbers which appear more thoroughly and establish a connection with the Lagrange and Markov spectra below~$3$.

\subsection{Markov triples}
\label{ssec:markovtriples}

We give a brief overview on Markov triples, highlighting the facts that will be used in the paper. 
For more details, see for example the book~\cite{Aig13}.
A \emph{Markov triple}\/ is a triple~$(a,b,c)$ of natural numbers (called \emph{Markov numbers}) 
satisfying the \emph{Markov equation},
\begin{equation}
\label{eq:markov}
a^2 + b^2 + c^2 \,=\, 3abc.
\end{equation}
We consider Markov triples to be ordered in the sense that we write~$(a,b,c)$ 
if~$a \geqslant b \geqslant c$. We denote the set of Markov triples by~$\mathfrak{M}$. 
Fixing two out of three numbers in a triple turns~\eqref{eq:markov} into a quadratic equation with two integer solutions. 
This means that, from one triple~$(a,b,c)$ one can obtain three new triples by 
so-called \emph{mutations}, 
\begin{equation}
	\label{eq:elmove}
	(a,b,c) \xrightarrow{\,c\;} (3ab-c,a,b), \quad
	(a,b,c) \xrightarrow{\,b\;} (3ac-b,a,c), \quad
	(a,b,c) \xrightarrow{\,a\;} (b,c,3bc-a).
\end{equation}
Therefore, Markov triples can be naturally arranged in a trivalent tree with root~$(1,1,1)$, and one can show that this tree contains all solutions. The tree is degenerate at the first two entries, in the sense that it contains repeated Markov triples (which we omit).

\begin{center}
\begin{tikzpicture}[level distance=1cm,
  level 3/.style={sibling distance=6cm},
  level 4/.style={sibling distance=3cm}]
  \node {$(1,1,1)$}
  child { node {$(2,1,1)$}
  child {node {$(5,2,1)$}
    child {node {$(13,5,1)$}
      child {node {$(34,13,1)$}}
      child {node {$(194,13,5)$}}
    }
    child {node {$(29,5,2)$}
    child {node {$(433,29,5)$}}
      child {node {$(169,29,2)$}}
    }}};
\end{tikzpicture}
\end{center}

The following conjecture was first stated by Frobenius in~1913 and is still open. 

\begin{conjecture}[Uniqueness Conjecture]
\label{conj:markovconjecture}
Every Markov triple is uniquely determined by its maximal Markov number.
\end{conjecture}

Note that if the uniqueness conjecture holds, then the statement of Theorem~\ref{thm:grwidth} can be simplified. Indeed, the relevant Markov triple~$(p,b,c) \in \fM$ is then simply the unique one 
in which~$p$ appears as the maximal entry. 
Let us check that the Markov triple we refer to in the statement of Theorem~\ref{thm:grwidth} 
is well-defined, independently of the Uniqueness conjecture. 
Let~$p \in \fm$ be a Markov number in a Markov triple~$\fm \in \fM$. 
By~$\bigwedge(\fm,p)$ we denote the set of all Markov triples which can be obtained from~$\fm$ 
by mutations preserving~$p$. This set forms a bivalent subtree of the Markov tree, 
since at each triple in~$\bigwedge(\fm,p)$ there are two mutations preserving~$p$. 
The subtree is rooted in a unique Markov triple containing~$p$ as maximal element -- 
this is the element we refer to in Theorem~\ref{thm:grwidth}. 
Indeed, take~$\fm' \in \bigwedge (\fm, p)$. If~$p$ is not the maximal entry in~$\fm$, 
we can repeatedly apply the third mutation in~\eqref{eq:elmove}. 
Under this mutation, the maximal entry strictly decreases and hence~$p$ will end up being maximal. 
The triple we obtain is unique, since the maximal entry strictly increases when applying one of the 
other two mutations in~\eqref{eq:elmove}. 

\subsection{Lagrangian pinwheels}

We consider the symplectic manifold~$(\CP^2,\omega_{\FS})$, where~$\omega_{\FS}$ is the Fubini--Study form normalized 
so that~$\int_{\CP^1} \omega_{\FS} = 1$. We are interested in \emph{Lagrangian pinwheels}\/ in~$\CP^2$. We give a brief overview here, see~\cite{ES18} and Section~\ref{ssec:atfs} for more details. 

\begin{definition}
Let~$p \in \NN_{\geqslant 2}$. The topological space obtained from the unit disk~$D$ 
by quotienting out the action of the group of the~$p$-th roots of unity on~$\pp D$ is called $p$-\emph{pinwheel}.
\end{definition}

For example, the~$2$-pinwheel is~$\RP^2$. The image of~$\pp D$ in the quotient is called the \emph{core circle}. For~$p > 2$ the~$p$-pinwheel is not smooth at points of the core circle. Pinwheels appear as skeleta of certain rational homology balls~$B_{p,q}$ such that the image of the interior of the disk is Lagrangian. The number~$q \in \mathbb{Z}$ describes the winding of the Lagrangian tangent plane about the core circle and is only determined up to sign and modulo~$p$. We call a pinwheel~$L_p \subset M$ \emph{Lagrangian $(p,q)$-pinwheel} if it admits a neighbourhood which is symplectomorphic to the normal form of~$L_p \subset B_{p,q}$ as discussed in~\cite[2.1]{ES18}. Evans--Smith~\cite{ES18} have studied Lagrangian embeddings of disjoint unions of pinwheels,
\begin{equation}
	\label{eq:pinwheelembedding}
	L_{p_1} \sqcup \dots \sqcup L_{p_N} 
	\hookrightarrow
	\CP^2.
\end{equation}

In fact, they give a complete classification of these embeddings in terms of the set of pairs~$\{(p_i,q_i)\}_{i \in 1,\ldots,N}$.

\begin{theorem}[Evans--Smith~\cite{ES18}]
\label{thm:evasmi}
Let there be a Lagrangian embedding of a disjoint union of pinwheels into~$\CP^2$ as in~\eqref{eq:pinwheelembedding}. Then~$N \leqslant 3$, and if~$N=3$, then~$(p_1,p_2,p_3)$ is a Markov triple. If~$N\in \{1,2\}$, then the~$p_i$ can be completed to a Markov triple.
\end{theorem}

See~\cite[Theorem 1.2]{ES18} for a more detailed statement.

Note that in Theorem~\ref{thm:grwidth} we restrict our attention to certain types of Lagrangian pinwheels, 
which we call \emph{Markov pinwheels}. This comes from the fact that it is unknown whether Lagrangian pinwheel embeddings are unique up to symplectomorphism
(except for~$L_2^{(2,1,1)} = \RP^2$), meaning that we do not know if Lagrangian pinwheels can be mapped to one another by an ambient symplectomorphism, see also Remark~\ref{rk:uniqueness}. 
Therefore, given a Markov triple~$\fm = (a,b,c) \in \fM$, we fix a particular embedding
\begin{equation}
\label{eq:markovembedding}
	L_a \sqcup L_b \sqcup L_c \hookrightarrow \CP^2.
\end{equation}
This embedding comes from almost toric fibrations of~$\CP^2$, see~\cite{Sym03,Via16,Via17}. 
Markov pinwheels are defined as certain subsets fibering over segments in the almost toric base diagrams of~$\CP^2$, 
see for example~\cite[Remark 2.9]{ES18}. 
This will be discussed in detail in Section~\ref{ssec:atfs}, and a precise definition of Markov pinwheels is given in 
Definition~\ref{def:markovpinwheel}. 
We denote Markov pinwheels obtained as connected components of the image of~$\eqref{eq:markovembedding}$ 
by~$L_a^{\fm}, L_b^{\fm}, L_c^{\fm} \subset \CP^2$. 

\begin{remark}[Uniqueness of Lagrangian pinwheels]
\label{rk:uniqueness}
{\rm
We prove in Lemma~\ref{lem:mutations} that the Markov pinwheels~$L_p^{\fm}$ and~$L_p^{\fm'}$ are equal as sets 
whenever~$\fm' \in \bigwedge(\fm,p)$. 
The related question whether~$L_p^{\fm}$ is symplectomorphic to every~$L_p^{\fm'}$ seems transcendent,
because it is equivalent to the Uniqueness Conjecture~\ref{conj:markovconjecture}. 
Indeed, the positive answer to this question would imply that~$p$ determines the extrinsic parameter~$q$ 
(up to sign and modulo~$p$), 
and it is known that this implies the Uniqueness Conjecture, 
see~\cite[Section 2.1]{ES18} and~\cite[Proposition~3.15]{Aig13}. 
A much more approachable conjecture is that Lagrangian~$(p,q)$-pinwheels in~$\CP^2$ 
are unique up to symplectomorphism. 
}
\end{remark}

\subsection{Main results}
Our main result is Theorem~\ref{thm:grwidth}. For the first few Markov triples we obtain,

\begin{center}
\renewcommand{\arraystretch}{1.6}
\begin{tabular}{ c||c|c|c|c|c  } 
 $(p,b,c)$ & $(2,1,1)$ & $(5,2,1)$ & $(13,5,1)$ & $(29,5,2)$ & $(433,29,5)$ \\ \hline
 $w_{\G} \bigl(\CP^2 \setminus L_p^{(p,b,c)} \bigr)$ & $ \frac{1}{2}$ & $\frac{2}{5}$ & $\frac{5}{13}$ & $\frac{10}{29}$ & $\frac{145}{433}$
\end{tabular}
\end{center}

\medskip
From Theorem~\ref{thm:grwidth}
it is not hard to see that~$\frac{1}{3} < w_{\G} \bigl(\CP^2 \setminus L_p^{(p,b,c)} \bigr) \leqslant \frac{1}{2}$ and that~$\frac{1}{3}$ 
is an accumulation point of the Gromov widths. 
In particular, every Markov pinwheel is a non-trivial Lagrangian barrier and every symplectic ball in~$\CP^2$ of capacity 
larger than~$\frac{1}{3}$ 
intersects infinitely many Lagrangian pinwheels. 
In Subsection~\ref{s:lattice} we discuss how the quantities~$\frac{bc}{a}$ are related to the base triangles of almost toric fibrations
of~$\CP^2$. 
Furthermore, in Section~\ref{s:numbers} we show that the accumulation points of the set of numbers obtained as relative Gromov widths are in bijection with the Lagrange and Markov spectra below~$3$.

\begin{corollary}
\label{cor:grwidth2}
Let~$\fm \in \fM$ be a Markov triple containing~$p_1,p_2$ with $p_1 > p_2$. Then
\begin{equation}
\label{eq:grwidth2}
	w_{\G} \bigl(\CP^2 \setminus (L_{p_1}^\fm \sqcup L_{p_2}^{\fm}) \bigr)
	\,=\,
	\frac{p_2 c}{p_1},
\end{equation} 
where~$(p_1,p_2,c) \in \mathfrak{M}$ is the unique Markov triple with $p_1>p_2$ and~$p_1>c$.
Let~$\fm = (p_1,p_2,p_3)$ be a Markov triple in decreasing order. Then
\begin{equation}
\label{eq:grwidth3}
	w_{\G} \bigl(\CP^2 \setminus (L_{p_1}^\fm \sqcup L_{p_2}^\fm \sqcup L_{p_3}^\fm) \bigr)
	\,=\,
	\frac{p_2 p_3}{p_1}.
\end{equation} 
\end{corollary}


Note that one inequality in Corollary~\ref{cor:grwidth2} is a direct consequence of Theorem~\ref{thm:grwidth}. 
The interest of the corollary lies the fact that the pinwheel associated to the largest Markov number determines the Gromov width of the complement, even when removing more than one pinwheel.


\medskip \noindent
{\bf Structure of the paper.}
The essence of the paper is contained in Subsection~\ref{ssec:atfs} and Section~\ref{s:proof}, which prove our main results. 
Subsection~\ref{ssec:atfs} contains background material from almost toric fibrations, 
and serves mainly to define and study Markov pinwheels. 
Section~\ref{s:proof} consists of the proof of Theorem~\ref{thm:grwidth} and Corollary~\ref{cor:grwidth2}. 
In Subsection~\ref{s:lattice} we give a geometric interpretation of the quantities~$w_{\G} \bigl(\CP^2 \setminus L_p^\fm \bigr) = \frac{bc}{p}$ in terms of the integral affine geometry of almost toric base diagrams. 
In Section~\ref{s:numbers} we recall the Lagrange and Markov spectra and show that the accumulation points of the set~$\{w_{\G}(\CP^2 \setminus L_p^{\fm})\}$ of relative Gromov widths are in bijection with these two spectra below~$3$.
In the appendix we argue that the immersed Lagrangians~$L_p^{\fm}$ can indeed be considered to be monotone.



%

\medskip \noindent
{\bf Acknowldegement.}
We started thinking about pinwheels as Lagrangian barriers after a suggestive remark
by Leonid Polterovich after the talk of JB at Tel Aviv in March 2022. 
Exciting discussions with 
Yasha Eliashberg at the~ITS of ETH~Z\"urich, with Jonny Evans, and with Grisha Mikhalkin
were key for the proof of Proposition~\ref{p:wG}.
We thank them all cordially.

\section{Almost toric fibrations}
\label{sec:atf}

Almost toric fibrations (ATFs) were introduced in~\cite{Sym03} and later used to great effect in~\cite{Via16,Via17}. We use the terminology and basic results from almost toric geometry without reference and would like to point out that our understanding of the topic was greatly influenced by~\cite{ES18},~\cite{Eva19} and~\cite{Eva21}, which we recommend for further details on ATFs. 
In~\cite{Via16}, Vianna introduced a family of almost toric fibrations of~$\CP^2$ which are in bijection with Markov triples. 
Furthermore, using the ATF-toolkit one can realize the algebraic mutations of Markov triples~\eqref{eq:elmove} geometrically 
in terms of the almost toric fibration. This means that there is a correspondence between the Markov tree and ATFs of~$\CP^2$,
under which each vertex corresponds to a different fibration and each edge corresponds to an ATF-mutation. 

\subsection{Markov pinwheels}
\label{ssec:atfs}

The main goal of this subsection is to define Markov pinwheels~$L^{\fm}_p \subset \CP^2$ in terms of almost toric fibrations. 
Let us stress again that it is unknown whether these pinwheels are unique up to Hamiltonian isotopy. 
We therefore first define them as sets, 
and then prove that the definition is invariant under mutations preserving~$p$ (Lemma~\ref{lem:mutations}) 
and that pinwheels can be made to lie arbitrarily close to the vertices of the ATF-base diagram by modifying the fibration structure 
(Lemma~\ref{lem:slides}). 
Both Lemmata will be used in the proof of Theorem~\ref{thm:grwidth}.

We fix almost toric fibrations
\begin{equation}
	\label{eq:pim}
	\pi_{\fm} \colon \CP^2 \rightarrow \Delta_{\fm}, \quad
	\fm \in \fM,
\end{equation}
where~$\pi_{\fm}$ is related to~$\pi_{\fm'}$ by an ATF-mutation if~$\fm$ and~$\fm'$ are connected by an edge in the Markov tree. 
Every base diagram~$\Delta_{\fm}$ is a triangle in the plane of area~$\frac{1}{2}$ (in our normalization of the symplectic form). 
Every triangle~$\Delta_{\fm}$ is decorated with marked points (crosses in the figures) 
-- which indicate a nodal fibre lying over that point, 
and dashed line segments connecting vertices of the triangle to the marked points -- which indicate the presence of monodromy. 
The dashed segments lie on lines which meet at a common point in the interior of~$\Delta_{\fm}$, 
which we call \emph{central point}\/ of~$\Delta_{\fm}$ and which lies at integral affine distance~$\frac{1}{3}$ to all edges. 
By \emph{nodal slides}\/ one can change the position of the marked points on the line spanned by the corresponding dashed segment. 
This means that when we fix~\eqref{eq:pim}
there are actually three additional nodal slide parameters~$\delta_1,\delta_2,\delta_3 > 0$,  
which measure the integral affine length of the dashed lines in the base diagram. 
For convenience we set~$\delta = \delta_1= \delta_2 = \delta_3$ and fix a small~$0 < \delta < \frac{1}{3}$ from now on. 
We denote the dashed line segments by~$\sigma_a(\delta),\sigma_b(\delta),\sigma_c(\delta) \subset \Delta_{\fm}$. 
By~\cite[Remark 2.9]{ES18} and~\cite[Section~2.1]{Eva19} there are pinwheels lying over the dashed lines 
in the base diagram\footnote{They are visible Lagrangians in the sense of~\cite[Section 7]{Sym03} and~\cite[Chapter 5]{Eva21}.}.

\begin{figure}[h] 
 \begin{center}
  \psfrag{c}{$(5,2,1)$} 
  \leavevmode\epsfbox{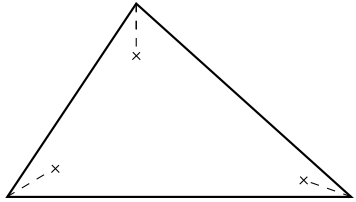}
 \end{center}
 \caption{$\sigma_a(\delta),\sigma_b(\delta),\sigma_c(\delta) \subset \Delta_{\fm}$}
 \label{fig:atd}
\end{figure}


\begin{definition}
\label{def:markovpinwheel}
Let~$\fm = (a,b,c) \in \fM$ be a Markov triple and let~$\pi_{\fm}$ be the almost toric fibration from~\eqref{eq:pim} with nodal slide parameter~$\delta$. Then we set the \emph{Markov pinwheels}~$L^{\fm}_a,L^{\fm}_b,L^{\fm}_c \subset \CP^2$ to be the Lagrangian pinwheels lying over~$\sigma_a(\delta),\sigma_b(\delta),\sigma_c(\delta) \subset \Delta_{\fm}$, respectively. 
\end{definition}

The Markov pinwheels~$L^{\fm}_a,L^{\fm}_b,L^{\fm}_c \subset \CP^2$ are, up to Hamiltonian isotopy, independent of the choice of the 
nodal slide parameter~$\delta$. Nodal slides induce Hamiltonian isotopies mapping~$\delta$-Markov pinwheels to~$\delta'$-Markov pinwheels. This follows from using the parallel transport maps coming from the local smoothing of orbifold singularities to 
rational homology balls~$B_{p,q}$ as described in~\cite{Eva19}. 
In particular,~\cite[Lemma 1.20]{Eva19} shows that pinwheels are mapped to pinwheels by the parallel transport. 
Since we do not use this fact, we do not give a detailed proof.

Recall that~$\bigwedge(\fm,p)$ is the set of Markov triples which can be obtained from~$\fm \in \fM$ by mutations preserving~$p \in \fm$. 

\begin{lemma}
\label{lem:mutations}
Let~$\fm' \in \bigwedge(\fm,p)$ be a Markov triple. Then the Markov pinwheels~$L^{\fm'}_p$ and~$L^{\fm}_p$ are equal as sets.
\end{lemma}

\proof
The claim follows straightforwardly from the ATF-toolkit. Let~$\fm = (p,b,c) \in \fM$ and suppose that we apply the mutation~$(p,b,c) \rightarrow \fm' = (3pc - b,p,c)$ which preserves~$p$. On the ATF-side, this corresponds to performing an ATF-mutation at the vertex corresponding to~$b$ (the Markov number which is eliminated in the triple). When applying an ATF-mutation to a vertex of~$\Delta_{\fm}$, the fibration structure changes only on an arbitrarily small neighbourhood of the segment defined by the dashed line at  the vertex in question and passing through the central point of~$\Delta_{\fm}$. In our normalization, the central point lies at integral affine distance~$\frac{1}{3}$ to all three edges of~$\Delta_{\fm}$. Therefore, since~$\delta < \frac{1}{3}$, the mutation at~$b$ does not change the fibration structure over~$\sigma_p(\delta)$. This means that~$L_p^{\fm} = L_p^{\fm'}$. Note that the same is true for mutations performed at~$c$, 
and hence the general claim follows from iterating.
\proofend

In the proof of Theorem~\ref{thm:grwidth}, we make use of the almost toric fibrations of~$\CP^2$. In particular, we need the fact that, after changing the fibration in a small neighbourhood of~$L_p^{\fm}$, we can assume that~$\pi_{\fm}(L_p^{\fm}) \subset \Delta_{\fm}$ lies arbitrarily close to the corresponding vertex.  As will become apparent in Section~\ref{s:proof}, using nodal slides is not sufficient for us, since the Hamiltonian isotopies they induce are not supported in a small neighbourhood of~$L_p^{\fm}$ but rather in the preimage of a small neighbourhood of~$\sigma_p(\delta)$.

\begin{lemma}
\label{lem:slides}
Let~$L_p^{\fm} \subset \CP^2$ be a Markov pinwheel and~$U$ be a neighbourhood thereof. Let~$\varepsilon > 0$. 
Then there exists an almost toric fibration~$\pi_{\fm}'$ which coincides with~$\pi_{\fm}$ outside of~$U$ and maps~$L_p^{\fm}$ to a segment~$\sigma \subset \Delta_{\fm}$ 
emanating from the corresponding vertex and of integral affine length smaller than~$\varepsilon$. 
\end{lemma}

Let us prepare the proof of this lemma by discussing a local model of the almost toric fibration~$\pi_{\fm}$ at its vertices.  
The pinwheel~$L_p^{\fm}$ has a neighbourhood symplectomorphic to a rational homology ball~$B_{p,q}$, 
where~$q$ is determined by the Markov triple~$\fm$. The space~$B_{p,q}$ can be equipped with an almost toric fibration, 
which yields an almost toric model for the fibration~$\pi_{\fm}$ at the vertex corresponding to~$p$. 
We refer to~\cite[Section~9.3]{Sym03} for a detailed discussion of the almost toric structure on~$B_{p,q}$ 
and to~\cite{Via16} for the local description in~$\CP^2$. 
In what follows, we heavily rely on the description of the spaces~$B_{p,q} = B_{1,p,q}$ discussed in~\cite{Eva19}. 
There is a conic fibration~$r \colon B_{p,q} \rightarrow \CC$ with a singular value at~$s_{\delta} = \sqrt{\delta/\pi} \in \CC$ 
and a special fibre at~$0 \in \CC$ such that the generic fibre is a~$p$-to-$1$ multiple cover of~$r^{-1}(0)$,
see~$(1)$ in \cite[Section~1.2.3]{Eva19}.
We chose~$s_{\delta}$ such that the segment~$\sigma_p(\delta)$ has integral affine length~$\delta$, see~\eqref{eq:g0} below. 
There is a Hamiltonian~$H$ on~$B_{p,q}$ which is compatible with~$r$ in the sense that its Hamiltonian flow preserves the fibres of~$r$. 
We normalize~$H$ such that all singular points of~$r$ lie in~$H^{-1}(0)$. 
The Markov pinwheel~$L_p^{\fm}$ is given as~$H^{-1}(0) \cap r^{-1}([0,s_{\delta}]) \subset B_{p,q}$. 
The connection between this conic fibration and the ATF-picture of~$B_{p,q}$ is described in~\cite[Example 2.11]{Eva19}. 
The coordinates of the almost toric projection of~$B_{p,q}$ are given by~$H$ and the function
\begin{equation}
\label{eq:g0}
G \colon B_{p,q} \rightarrow \RR, \quad
x \mapsto g(r(x)), \quad
g(z) = \pi \vert z \vert^2
\end{equation}
defined on the base of the conic fibration. 
The image of the map~$(H,G) \colon B_{p,q} \rightarrow \R^2$ is the upper half-plane, which looks very different from the vertex 
in question in~$\Delta_{\fm}$; however, the two pictures are related by a \emph{mutation}, 
see~\cite[Section 2.1.4]{Eva19}. 
In order to prove Lemma~\ref{lem:slides}, we slightly modify this fibration structure. 

\medskip \noindent
\emph{Proof of Lemma~\ref{lem:slides}: }
Let~$U$ be a neighbourhood of the pinwheel and~$\varepsilon > 0$. 
By~$W_{\lambda} \subset \CC$ we denote all points at distance smaller than~$\lambda > 0$ from 
the segment~$[0,s_{\delta}] \subset \RR \subset \CC$. 
Fix~$\lambda$ small enough such that the set~$V_{\lambda}= H^{-1}(-\lambda,\lambda) \cap r^{-1}(W_{\lambda})$ is contained in~$U$. 
Since~$L^{\fm}_p = H^{-1}(0) \cap r^{-1}([0,s_{\delta}])$, the set~$V_\lambda$ is a neighbourhood of~$L_p^{\fm}$. 
We modify the fibration inside~$V_{\lambda} \subset U$, 
and in terms of the coordinate functions~$(H,G)$ only the~$G$-coordinate will be modified.

\begin{figure}[h] 
 \begin{center}
  \psfrag{W}{$W_\lambda$} \psfrag{0}{$[0,s_\delta]$} 
	\psfrag{L1}{$\mbox{Level sets of $g$}$} 	\psfrag{L2}{$\mbox{Level sets of $g \circ \gf_1$}$} 
  \leavevmode\epsfbox{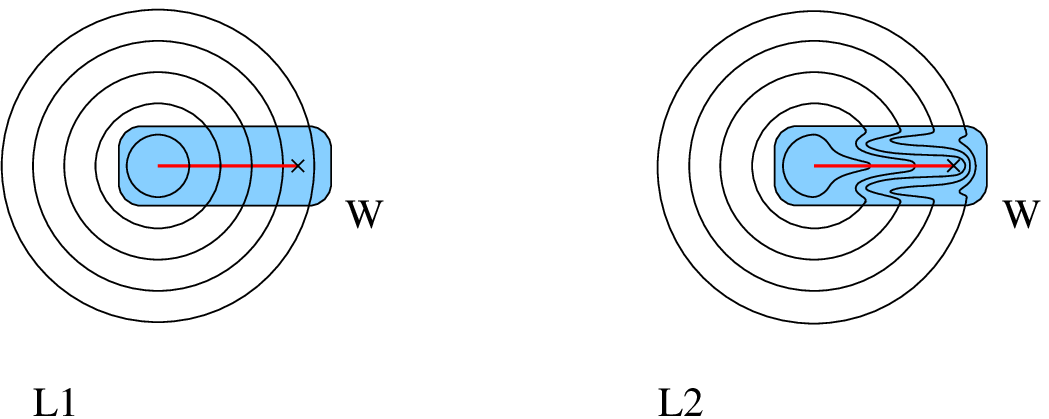}
 \end{center}
 \caption{}
 \label{fig:levelsets}
\end{figure}

%

For this, we choose a Hamiltonian isotopy~$\{\varphi_t\}_{t \in [0,1]}$ of the plane~$\CC$ with the following properties,
\begin{enumerate}
\item $\supp(\varphi_t) \subset W_{\lambda}$ \vspace{0.2em}
\item $0 \leqslant (g \circ \varphi_1)\vert_{[0,s_{\delta}]} < \varepsilon $.
\end{enumerate}
Recall from~\eqref{eq:g0} that~$g(z) = \pi \vert z \vert^2$, and see Figure~\ref{fig:levelsets} for a sketch of such a Hamiltonian isotopy. 
Let~$b \colon \R \rightarrow [0,1]$ be a smooth bump function with~$\supp b \subset (-\lambda,\lambda)$ and~$b(0) = 1$. 
Set
\begin{equation}
	\widetilde G(x) \,=\, (g \circ \varphi_{b(H(x))} \circ r)(x).
\end{equation}
Note that~$\widetilde G$ coincides with~$G$ outside of~$V_{\lambda}$. 
Furthermore, the new function~$\widetilde G$ still defines a Hamiltonian~$S^1$-action and commutes with~$H$. 
Indeed, the value~$\widetilde G(x)$ depends only on~$H(x)$ and~$r(x)$, both of which are preserved by 
the Hamiltonian flow generated by~$H$, and therefore~$d\widetilde G(X_H) = \{\widetilde G,H\} = 0$. 
Since this modification of~$G$ to~$\widetilde G$ is local in a neighbourhood of~$L^{\fm}_p$, 
it can be carried out such that the given fibration~$\pi_{\fm}$ is altered to the desired fibration~$\pi_{\fm}'$. 
The fact that the segment~$\sigma = \pi_{\fm}'(L_p^{\fm}) \subset \Delta_{\fm}$ has affine length~$< \varepsilon$ 
follows from condition~$(2)$ imposed on~$\varphi_1$ together with the fact that~$L_p^{\fm} = H^{-1}(0) \cap r^{-1}([0,s_{\delta}])$ 
in the local model.
\proofend

\subsection{Geometric interpretation of the relative Gromov widths}
\label{s:lattice}
The quantity~$\frac{bc}{a}$ for~$\fm = (a,b,c) \in \mathfrak{M}$ appearing in Theorem~\ref{thm:grwidth} has a geometric interpretation as the \emph{lattice width} of the integral affine geometry of the corresponding almost toric base triangle~$\Delta_{\fm}$. 

\begin{definition}
The \emph{lattice width}\/ of a compact subset~$A \subset \RR^n$ is defined as 
\begin{equation}
	\label{eq:latticewidth}
	w_{\aff}(A) 
	\,=\, 
	\min_{\xi \in \ZZ^n  \setminus \{0\}} 
	\left( \max_{x,x' \in A} \langle x' - x , \xi \rangle \right).
\end{equation}
\end{definition}

\begin{proposition}
\label{prop:affwidth}
Let~$\fm = (a,b,c) \in \fM$ be an ordered Markov triple. Then 
\begin{equation}
	w_{\aff}(\Delta_{\fm}) \,=\, \frac{bc}{a}.
\end{equation}
\end{proposition}

\proof
Let~$E$ be the longest edge of the triangle~$\Delta_{\fm}$ and denote by~$\ell$ its integral affine length. Up to applying a lattice automorphism and a translation, we may assume that~$E$ lies on 
the coordinate line spanned by~$(1,0)$ with one vertex at~$(0,0)$ and one vertex at~$(\ell,0)$. 
We claim that~$w_{\aff}(\Delta_{\fm})$ is realized by~$\xi = (0,1)$, 
meaning that it is equal to the height~$h$ of the triangle~$\Delta_{\fm}$ in this normal form. 

Let us first show that this height equals~$\frac{bc}{a}$. 
We have~$2 \area \Delta_{\fm} = h \ell$ and $\area \Delta_{\fm} = \frac{1}{2}$, hence $h \ell = 1$.
There is a constant $\lambda >0$ such that the length of the edge opposite to the vertex of weigth $a,b,c$
is $\lambda a^2, \lambda b^2, \lambda c^2$, respectively, see~\cite[Section 2]{Via16}.
In view of our ordering~$a \geqslant b \geqslant c$ we have
$\ell = \lambda a^2$. 
Since the perimeter of~$\Delta_{\fm}$ is~$3$, we can use the Markov equation to find that 
$\lambda = \frac{1}{abc}$, whence~$\ell = \frac{a}{bc}$, and so $h= \frac 1 \ell = \frac{bc}{a}$. 

Let us now prove that~$\xi = (0,1)$ realizes the minimum in~\eqref{eq:latticewidth}. 
Assume that~$\xi = (\xi_1,\xi_2)$ with~$\xi_1 \neq 0$. 
Since 
$\displaystyle \max_{x,x' \in \Delta_{\fm}} \langle x' - x , \xi \rangle = \max_{x,x' \in \Delta_{\fm}} \langle x' - x , - \xi \rangle$, 
we can assume that $\xi_1 >0$.
We choose the points~$x_0 = (0,0)$ and~$x_0' = (\ell,0)$ in~$\Delta_{\fm}$ and note that 
\begin{equation}
	\max_{x,x' \in \Delta_{\fm}} \langle x' - x , \xi \rangle
	\,\geqslant\, \langle x_0' - x_0 , \xi \rangle
	\,=\, \langle (\ell,0) , \xi \rangle
	\,=\, \ell \xi_1 \,\geqslant\, \ell \,\geqslant\, h,
\end{equation}
where the last inequality follows from the fact that~$\ell \geqslant 1$
(because $\ell$ is the length of the longest edge of the triangle~$\Delta_{\fm}$ of perimeter~$3$)
together with~$h \ell = 1$.  
Thus~$\xi = (0,1)$ realizes the minimum.
\proofend

In the toric case, the lattice width of the Delzant polytope is intimately related to the Gromov width of the corresponding space.

\begin{theorem}[Chaidez--Wormleighton~\cite{ChWo20}]
\label{thm:widthinequality}
Let~$(M,\omega)$ be a toric symplectic manifold with Delzant polytope~$\Delta = \mu (M)$. 
Then
\begin{equation}
	\label{eq:widthinequality}
	w_{\G}(M) \,\leqslant\, w_{\aff}(\Delta).
\end{equation}
\end{theorem}

This was conjectured earlier in~\cite{AHN19} and \cite{HLS21}. 
While the upper bound in~\eqref{eq:widthinequality} is not always optimal, it is optimal in many examples. 
Proposition~\ref{prop:affwidth} shows that Theorem~\ref{thm:widthinequality} fails in the almost toric setting. 
Indeed, by Proposition~\ref{prop:affwidth} and the Markov equation  
we have~$w_{\aff}(\Delta_\fm) < 1$ for all~$\fm \neq (1,1,1)$, but~$w_{\G}(\CP^2)=1$. 
Theorem~\ref{thm:grwidth} suggests that 
the lattice width 
of almost toric bases 
should be interpreted as 
the \emph{relative} Gromov width 
for suitable choices of a Lagrangian.

\section{Relation with the Lagrange spectrum, and order in the Markov tree} 
\label{s:numbers}

\subsection{Accumulation at the Lagrange spectrum below 3}

Dirichlet proved that for every irrational number $\alpha$
there exist infinitely many rational numbers $\frac pq$ such that
$$
\left| \alpha - \frac pq \right| \,\leq\, \frac{1}{q^2} ,
$$
and Roth proved that, in general, one cannot improve on the exponent~$2$.
So one tries to improve the constant~$1$:
The Lagrange number of~$\alpha$ is defined as the supremum over $L \geq 1$
such that
$$
\left| \alpha - \frac pq \right| \,\leq\, \frac{1}{L \2 q^2} 
$$
for infinitely many $\frac pq$.
The set $\cl = \left\{ L(\alpha) \mid \mbox{$\alpha$ irrational} \right\} \subset [1,\infty]$
is called the Lagrange spectrum. 
Hurwitz showed that for the Golden ratio, $L \bigl( (1+\sqrt{5})/2 \bigr) = \sqrt 5$, 
and that this is the minimum of~$\cl$.
The next numbers in~$\cl$ are $\sqrt 8$ and $\sqrt{221}/5$.

Let $f(x,y) = a \2 x^2 + b \2 xy + c \2 y^2$ with $a,b,c \in \RR$
be a quadradic form on~$\RR^2$, and assume that $f$ is indefinite:
$\Delta (f) = b^2-4ac >0$. 
Set 
$$
m(f) \,=\, \inf \left\{ |f(x,y)| : f(x,y) \neq 0, \; (x,y) \in \ZZ^2 \right\} ,
$$
and define the Markov spectrum as the set of numbers $\sqrt{\Delta (f)} / m(f)$
obtained in this way.

Markov proved that the Markov spectrum below~$3$ is equal to the Lagrange spectrum below~$3$. 
Furthermore, he showed that this set is parametrized by the set of Markov numbers~$\cm$ 
(as defined in Subsection~\ref{ssec:markovtriples}) by the bijection
$$
\lambda \colon \cm \to \cl \cap [1,3), \quad \lambda (a) = \sqrt{9-4/a^2} \, .
$$ 
We refer to~\cite{Aig13}, \cite{Bo07}, and \cite[Chapter 2]{Ca57} for much more on this beautiful story.

Now fix a Markov number $a$.
Take a Markov triple~$\fm$ with $a \in \fm$, 
and as in the introduction consider the subtree $\bigwedge (\fm,a)$ of Markov triples
obtained from~$\fm$ by mutations that preserve~$a$.
Hence $\bigwedge (\fm, 1)$ is the left-most branch of the Markov tree 
(the Fibonacci branch, formed by the triples $(F_{2n+1}, F_{2n-1}, 1)$ with odd-index Fibonacchi numbers),
$\bigwedge (\fm,2)$ is the right-most branch of the Markov tree 
(the Pell branch, formed by the triples $(P_{2n+1}, P_{2n-1}, 2)$ with odd-index Pell numbers),
and for $a \geq 5$ the tree $\bigwedge (\fm,a)$ is a $\bigwedge$-shaped bivalent subtree of the Markov tree,
cf.\ Figure~\ref{fig:order5}.
We recall that $\bigwedge (\fm,a)$ is independent of~$\fm$ exactly if the Markov uniqueness conjecture holds true for~$a$.

For $\fm = (p,b,c) \in \fM$ with $p \geq b,c$ we abbreviate the capacity
\begin{equation} \label{e:wfm}
w (\fm) \,=\, w_{\G} (\CP^2 \setminus L^\fm_{p}) .
\end{equation}
Furthermore, we say that a sequence of triples $(\fm_n)$ in $\bigwedge (\fm,a)$ is decreasing if
the distance of $\fm_n$ to~$(1,1,1)$ in the Markov tree is non-decreasing.

\begin{proposition} \label{p:limit}
For any decreasing sequence $(\fm_n)$ in~$\bigwedge (\fm,a)$,
\begin{equation} \label{e:limit}
\lim_{n \to \infty} w (\fm_n) \,=\, \frac{2}{3+\lambda(a)} \,=\, \frac{2}{3+ \sqrt{9-4/a^2}} .
\end{equation}
\end{proposition}

In words: if we go down along a subtree $\bigwedge (\fm,a)$,
then the Gromov width~$w (\fm)$ of the complement of the pinwheel associated to the triple~$\fm$
converges to the inverse of the arithmetic mean of~$3$ and the Lagrange number corresponding to~$a$ 
(that increases to~$3$ for $a$ large).

Note that the limit in \eqref{e:limit} does not depend on the choice of the triple~$\fm$ containg~$a$.
We therefore obtain a bijection
$$
\cm \to \left\{ \frac{2}{3+\lambda(a)} \;\bigg|\; a \in \cm \right\} , \qquad a \mapsto \frac{2}{3+\lambda (a)}
$$
that takes the set of Markov numbers to the set of accumulation points of $\left\{ w(\fm) \mid \fm \in \cm \right\}$.

\medskip
\noindent
{\it Proof of Proposition~\ref{p:limit}:}
It suffices to prove the claim for the case that $(\fm_n)$ consists either of the Fibonacci branch, 
or the Pell branch, or the left or the right branch of $\bigwedge (\fm, a)$.
After omitting $\fm_1$ and $\fm_2$, we have that $\fm_n$ is of the form $(c_n,b_n,a)$
with $c_n > b_n > a$.
Using Theorem~\ref{thm:grwidth} and the Markov equation $c_n^2 - (3 a b_n)\2 c_n + (a^2+b_n^2) =0$,
we compute
$$
w(\fm_n) \,=\, \frac{a \2 b_n}{c_n} \,=\, \frac{2}{3+\sqrt{9-4/a^2-4/b_n^2}} .
$$
Therefore, 
$$
\lim_{n \to \infty} w(\fm_n) \,=\, \frac{2}{3+\sqrt{9-4/a^2}} \,=\, 
\frac{2}{3+\lambda (a)},
$$
as claimed.
\proofend		

\subsection{Order of $\{ w(\fm)\}$ in the Markov tree}
	
Order the set $\{ w(\fm)\}$ defined by~\eqref{e:wfm} as a decreasing sequence.
How is this sequence arranged in the Markov tree?

Given an ordered Markov triple $\fm = (a,b,c)$ we now write $\bigwedge (\fm)$
for the bivalent tree $\bigwedge (\fm, a)$ with apex~$\fm$.
The following proposition makes Proposition~\ref{p:limit} more precise.
						
\begin{figure}[h] 
 \begin{center}
  \psfrag{c}{$(5,2,1)$} 
  \psfrag{l1}{$(13,5,1)$}  \psfrag{r1}{$(29,5,2)$} 
	\psfrag{l2}{$(194,13,5)$}  \psfrag{r2}{$(433,29,5)$} 
  \psfrag{l3}{$(2897,194,5)$}  \psfrag{r3}{$(6466,433,5)$} 
  \leavevmode\epsfbox{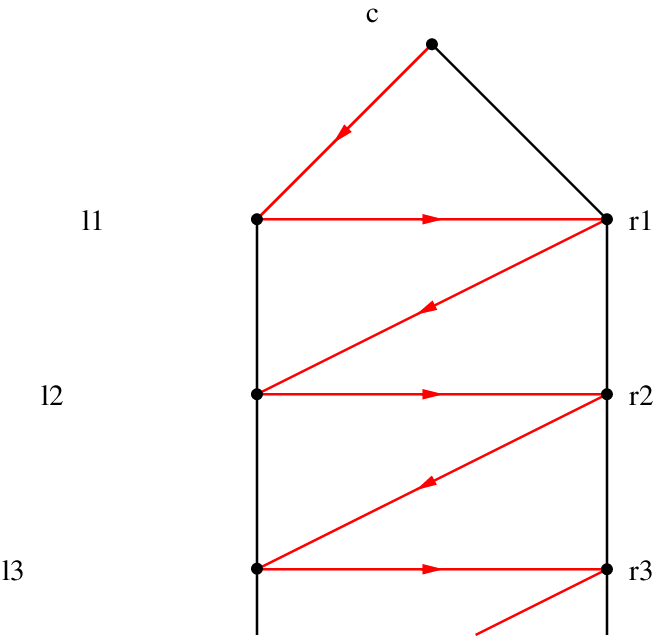}
 \end{center}
 \caption{The order of $w(\fm)$ along $\bigwedge \bigl( (5,2,1) \bigr)$}
 \label{fig:order5}
\end{figure}

\begin{proposition} \label{p:order}						
{\rm (i)}
Along the Fibonacci branch and along the Pell branch, the sequence formed by the~$w (\fm)$ 
is strictly decreasing. 						

\m
{\rm (ii)}
On each subtree $\bigwedge (\fm)$, $a \geq 5$,
the $w (\fm)$ strictly decrease if we go alternating between the two branches, 
going along increasing maximal Markov numbers:
start at the apex~$\fm$, go down left, go right, go down left, go right, etc.
See Figure~\ref{fig:order5} for~$a=5$.
\end{proposition}

%

\proof
Let $\fm = (a,b,c)$ with $a \geq b \geq c$ be a Markov triple.
As in the previous proof we deduce from Theorem~\ref{thm:grwidth} and from the Markov equation
that
\begin{equation} \label{e:wa}
w (\fm) \,=\, \frac{bc}{a} \,=\, \frac{2}{3+\sqrt{9-4/c^2-4/b^2}} .
\end{equation}

\smallskip
(i)
For the Fibonacci branch, $c=1$ and $a=F_{2n+1}$, $b=F_{2n-1}$, whence
$$
w \bigl((F_{2n+1}, F_{2n-1},1)\bigr) \,=\, \frac{2}{3+ \sqrt{5-4/F_{2n-1}^2}} .
$$
Similarly, for the Pell branch, $c=2$ and $a=P_{2n+1}$, $b=P_{2n-1}$, whence
$$
w \bigl((P_{2n+1}, P_{2n-1}, 2)\bigr) \,=\, \frac{2}{3+ \sqrt{8-4/P_{2n-1}^2}} .
$$
Since the sequences $(F_{2n+1})$ and $(P_{2n+1})$ are strictly increasing, 
assertion~(i) of the proposition follows.

\smallskip
(ii)
Fix an ordered Markov triple $\fm = (a,b,c)$ with $a \geq 5$. Then $a>b>c$.
The inequality that we have to check between the capacities of the apex and the largest descendent obtained by mutation is 
$$
\frac{bc}{a} \,>\,\frac{ac}{3ac-b} ;
$$
it follows from the Markov equation for $(a,b,c)$.
To save space, we draw $\bigwedge (\fm)$ deprived from the apex horizontally. 
The beginning is
\begin{equation*} 
\xymatrixcolsep{3pc}
\xymatrix{ 
(f_1,a,b) \ar@{-}^{b}[r]  \ar@{-->}[dr]  & (f_2,f_1,a) \ar@{-}^{f_1 \quad \;\; }[r] \ar@{-->}[dr]
& (3af_2-f_1, f_2,a) 
\\
(g_1,a,c) \ar@{-->}[u] \ar@{-}^{c}[r]  & (g_2,g_1,a) \ar@{-->}[u] \ar@{-}^{g_1 \quad \;\; }[r]& (3ag_2-g_1, g_2,a) 
}
\end{equation*}
where we abbreviated $f_1 = f_1(b,c) = 3ab-c$, $f_2 = f_2(b,c) = 3af_1-b$,
and $g_j(b,c) = f_j (c,b)$.
We must show that the associated capacities are ordered as indicated by the dashed
arrows:
$$
\frac{ac}{g_1} \,>\, \frac{ab}{f_1} \,>\, \frac{g_1a}{g_2} 
\,>\, \frac{f_1a}{f_2} \,>\, \frac{g_2a}{3ag_2-g_1}. 
$$
The first and third~$>$ are equivalent to~$b>c$, 
the second~$>$ is equivalent to $2b < 3ac$, and the fourth~$>$
is equivalent to $f_1/b < g_2/g_1$ and hence to the second~$>$.

The sequel of $\bigwedge (\fm)$ looks like 
\begin{equation*} 
\xymatrixcolsep{3pc}
\xymatrix{ 
\ar@{-}[r] (f_n,f_{n-1},a)  \ar@{-->}[dr]  & (3af_n-f_{n-1},f_n,a) \ar@{-}[r] \ar@{-->}[dr]
& \bigl(3a (3af_n-f_{n-1}) -f_n, 3af_n-f_{n-1},a \bigr) 
\\
\ar@{-}[r] (g_n,g_{n-1},a) \ar@{-->}[u] \ar@{-}[r]  & (3ag_n-g_{n-1},g_n,a) \ar@{-->}[u] \ar@{-}[r]&  
\bigl( 3a (3ag_n-g_{n-1}) -g_n, 3ag_n-g_{n-1},a \bigr)
}
\end{equation*}
where $f_{n-1} = f_{n-1}(b,c)$ and $f_{n} = f_n(b,c)$ are certain polynomials in~$b,c$
and $g_j = g_j(b,c) = f_j(c,b)$,
and where all Markov moves are done at the middle entry.

By induction, we assume that the capacity inequalities at the first two dashed arrows
hold true:
\begin{eqnarray} \label{e:ind}
\frac{f_{n-1}a}{f_n} \,<\, \frac{g_{n-1}a}{g_n} 
\quad \mbox{ and } \quad
\frac{g_n a}{3ag_n-g_{n-1}} \,<\, \frac{f_{n-1}a}{f_n} . 
\end{eqnarray}
We need to show that
$$
\frac{f_n a}{3af_n-f_{n-1}} \,<\, \frac{g_na}{3ag_n - g_{n-1}}
\quad \mbox{ and } \quad
\frac{(3a g_n - g_{n-1}) a}{3a (3a g_n-g_{n-1}) -g_n} \,<\, \frac{f_n a}{3a f_n - f_{n-1}} . 
$$
The left (resp.\ right) of these inequalities transforms to the left
(resp.\ right) inequality in~\eqref{e:ind}. 
\proofend

Define the essential part $\bigwedge_{\ess} (\fm)$ as the subset of~$\bigwedge (\fm) = \bigwedge (\fm,a)$
consisting of the triples with minimum~$a$.
The Markov tree is then the disjoint union of the~$\bigwedge_{\ess} (\fm)$.
For the Fibonacci branch, $\bigwedge (1,1,1) = \bigwedge_{\ess} (1,1,1)$,
for the Pell branch, $\bigwedge_{\ess} (2,1,1)$ is obtained from $\bigwedge (2,1,1)$
by removing $(2,1,1)$ and~$(5,2,1)$,
and $\bigwedge_{\ess} (\fm)$ with $a \geq 5$ is obtained from $\bigwedge (\fm)$
by removing the apex and its two children.
Restricting the decreasing sequences from Proposition~\ref{p:order}
we obtain decreasing sequences $w \bigl( \bigwedge_{\ess} (\fm)\bigr)$.
In view of Proposition~\ref{p:order} (and the Uniqueness conjecture) 
one can hope that the order of the decreasing sequence formed
by all the~$w (\fm)$ is simply given by juxtaposing the sequences $w \bigl( \bigwedge_{\ess} (\fm) \bigr)$	
in the order of the maximal elements of the~$\fm$.
As we shall see, 
this is indeed so for the first $32$~triples, but then fails for the first time.

Let $(m_1, m_2, m_3, \dots) = (1,2,5, \dots)$ be the increasing sequence of
Markov numbers.  
The Uniqueness conjecture is known to hold for $n \leq 18 \2 906$, i.e., 
for Markov numbers below~$10^{140}$,
see~\cite[Section 4]{Ba96}.
Since we will look only at the first $500$~Markov numbers, we can therefore write
$\bigwedge_{\ess} (m_n)$ instead of $\bigwedge_{\ess} (\fm)$, 
where $m_n$ is the maximum of the triple~$\fm$.
By~\eqref{e:wa}, the capacities in $w \bigl( \bigwedge_{\ess} (m_n) \bigr)$
form a sequence that strictly decreases to 
$$
\frac{2}{3+\sqrt{9-4/m_n^2}} .
$$
For $n'>n$, the largest capacity in 
$w \bigl( \bigwedge_{\ess} (m_{n'}) \bigr)$
is 
$$
\frac{2}{3+\sqrt{9-4/m_n^2 - 4/b_{n'}^2}} ,
$$
where $b_{n'}$ is the second smallest (after $m_{n'}$) Markov number
appearing in a triple of~$\bigwedge_{\ess} (m_{n'})$.
For instance, $(m_3,b_3) = (5,13)$ and $(m_4, b_4) = (13,34)$.
In the decreasing sequence formed by the $w(\fm)$, the sequence
$w \bigl( \bigwedge_{\ess} (m_{n}) \bigr)$ will therefore appear entirely before
any capacity from $w \bigl( \bigwedge_{\ess} (m_{n'}) \bigr)$ if and only if
\begin{equation} \label{e:nn'}
\frac{1}{m_n^2} \,\geq\, \frac{1}{m_{n'}^2} + \frac{1}{b_{n'}^2} .
\end{equation}
Our further discussion is based on little computer codes, for instance one 
associating to the Markov number~$a$ the first Markov triple in which $a$
appears, and on the lists 
\begin{eqnarray*}
\texttt{ oeis.org/A002559/b002559.txt } \\
\texttt{ oeis.org/A000045/b000045.txt } \\
\texttt{ oeis.org/A000129/b000129.txt }
\end{eqnarray*}
giving the first 1000 Markov, Fibonacchi, and Pell numbers.
Things start well: \eqref{e:nn'} holds for all $n'>n$
whenever $n \leq 32$. We therefore have:

\begin{proposition}
The beginning of the decreasing sequence formed by all the~$w (\fm)$ is obtained by juxtaposing
the sequences $w \bigl( \bigwedge_{\ess} (m_n) \bigr)$, $n \leq 32$:
$$
\textstyle
w \bigl( \bigwedge (1) \bigr), \;  w \bigl( \bigwedge_{\ess} (2) \bigr), \; w \bigl( \bigwedge_{\ess} (5) \bigr), \dots , 
 \; w \bigl( \bigwedge_{\ess} (m_{32}) \bigr) . 
$$
\end{proposition}

For $n=33$, however, \eqref{e:nn'} is violated for the first time,
because $m_{33}$ and~$m_{34}$ are very close:
We have $m_{33} = 195 \, 025$ and $m_{34} = 196 \, 418$, 
and the middle Markov numbers in their triples are 
$b_{33} = 1 \, 136 \, 689$ and $b_{34} = 514 \, 229$.
Hence the largest capacity $w_1 (m_{34})$ in $w \bigl( \bigwedge_{\ess} (m_{34}) \bigr)$
is larger than all capacities in $w \bigl( \bigwedge_{\ess} (m_{33}) \bigr)$.
All other capacities in $w \bigl( \bigwedge_{\ess} (m_{34}) \bigr)$ are smaller than
all capacities in $w \bigl( \bigwedge_{\ess} (m_{33}) \bigr)$.
The expected order 
$$
\textstyle
w \bigl( \bigwedge_{\ess} (m_{33}) \bigr), \; w \bigl( \bigwedge_{\ess} (m_{34}) \bigr)
$$
must therefore be replaced by  
$$
\textstyle
w_1 \bigl( \bigwedge_{\ess} (m_{34}), \; w \bigl( \bigwedge_{\ess} (m_{33}) \bigr), \;
w_{\geq 2} \bigl( \bigwedge_{\ess} (m_{34}) . 
$$
Note that $m_{33} = P_{15}$ is a Pell number and that $m_{34} = F_{27}$
is a Fibonacci number.

\begin{figure}[h] 
 \begin{center}
  \psfrag{w32}{$w \bigl( \bigwedge_{\ess} (m_{32}) \bigr)$} 
  \psfrag{w33}{$w \bigl( \bigwedge_{\ess} (m_{33}) \bigr)$} 
	\psfrag{w34}{$w \bigl( \bigwedge_{\ess} (m_{34}) \bigr)$} 
	\psfrag{R}{$\RR$} \psfrag{13}{$\frac 13$}    
  \leavevmode\epsfbox{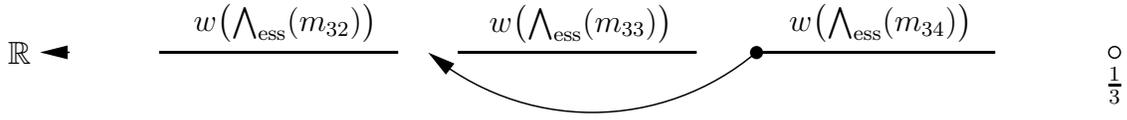}
 \end{center}
 \caption{The ``accident" at $n=33$}
 \label{fig:n33}
\end{figure}

%

Up until $n=450$, there are 16 further irregularities, namely at 
$(m_{n} , m_{n+1})$ with
$$
n = 37, \; 42, \; 
104, \; 112, \; 118, \; 120, \; 214, \; 227, \; 309, \; 353, \; 382, \; 
   400, \;416, \; 450
$$
and at the two places $(m_{n}, m_{n+2})$ with
$$
n= 369, \; 433 .
$$
The irregularities at $(m_{n}, m_{n+1})$ are exactly as for $(m_{33}, m_{34})$:
The expected order must be replaced by 
$$
\textstyle
w_1 \bigl( \bigwedge_{\ess} (m_{n+1}) \bigr), \; w \bigl( \bigwedge_{\ess} (m_{n}) \bigr), \;
w_{\geq 2} \bigl( \bigwedge_{\ess} (m_{n+1})\bigr) ,
$$
i.e., only the first capacity of the higher essential sequence must be swapped with the previous sequence.
Similarly, for $n= 369, \, 433$ the expected order
$$
\textstyle
w \bigl( \bigwedge_{\ess} (m_{n}) \bigr), \; w \bigl( \bigwedge_{\ess} (m_{n+1}) \bigr),
\; w \bigl( \bigwedge_{\ess} (m_{n+2}) \bigr)
$$
must be replaced by 
$$
\textstyle
w_1 \bigl( \bigwedge_{\ess} (m_{n+2}) \bigr), \; w \bigl( \bigwedge_{\ess} (m_{n}) \bigr), \;
w \bigl( \bigwedge_{\ess} (m_{n+1}) \bigr), \;
w_{\geq 2} \bigl( \bigwedge_{\ess} (m_{n+2}) \bigr),
$$
i.e., only the first capacity of the two higher essential sequence must be swapped with the previous two sequences.
We do not know whether for $n>450$ there are other kinds of irregularities.

\begin{figure}[h] 
 \begin{center}
  \psfrag{w32}{$w \bigl( \bigwedge_{\ess} (m_{n}) \bigr)$} 
  \psfrag{w33}{$w \bigl( \bigwedge_{\ess} (m_{n+1}) \bigr)$} 
	\psfrag{w34}{$w \bigl( \bigwedge_{\ess} (m_{n+2}) \bigr)$} 
  \psfrag{R}{$\RR$} \psfrag{13}{$\frac 13$}    
  \leavevmode\epsfbox{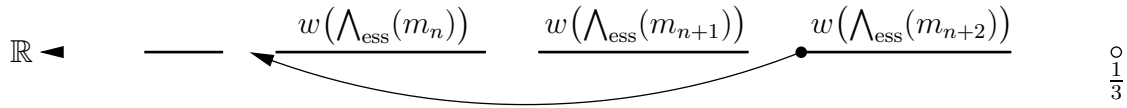}
 \end{center}
 \caption{The accidents at $n= 369, \, 433$}
 \label{fig:n369}
\end{figure}

%

The irregularities seem to come from Markov numbers in ``small branches":
All of the right numbers $m_{n+1}$ and $m_{n+2}$ (except $m_{119}$ and $m_{310}$)
are Fibonacci numbers, 
and most of the~$m_{n}$ are Pell numbers or very close to the Pell branch.

The above analysis gives 
a complete description of the decreasing sequence of those~$w (\fm)$ 
with $w (\fm) - 1/3 \geq 2 \cdot 10^{-44}$.
\proofend


\section{Proof of Theorem~\ref{thm:grwidth} and of Corollary~\ref{cor:grwidth2}}
\label{s:proof}

We shall show that for every Markov triple $\fm = (a,b,c)$ with $a \geqslant b \geqslant c$,
$$
w_{\G} \left( \CP^2 \setminus L_{a}^{\fm} \right) \,\leq\,
\frac{bc}{a} \,\leq\, 
w_{\G} \left( \CP^2 \setminus (L_{a}^\fm \sqcup L_{b}^\fm \sqcup L_{c}^{\fm}) \right) .
$$
Theorem~\ref{thm:grwidth} and Corollary~\ref{cor:grwidth2} then follow from the monotonicity of the Gromov width.

As in Section~\ref{ssec:atfs}
we consider the almost toric fibration $\pi_\fm \colon \CP^2 \to \Delta_\fm$.

\m \noindent
{\bf Step 1.}
$w_{\G} \bigl( \CP^2 \setminus (L_{a}^\fm \sqcup L_{b}^\fm \sqcup L_{c}^{\fm}) \bigr) \geqslant \frac{bc}a$

\begin{figure}[h] 
 \begin{center}
  \psfrag{1}{$x_1$}   \psfrag{2}{$x_2$} 
  \psfrag{a}{$\frac{a}{bc}$}  \psfrag{b}{$\frac{bc}{a}$} 
	\psfrag{M}{$M_k$} 	\psfrag{T}{$T$}  
  \leavevmode\epsfbox{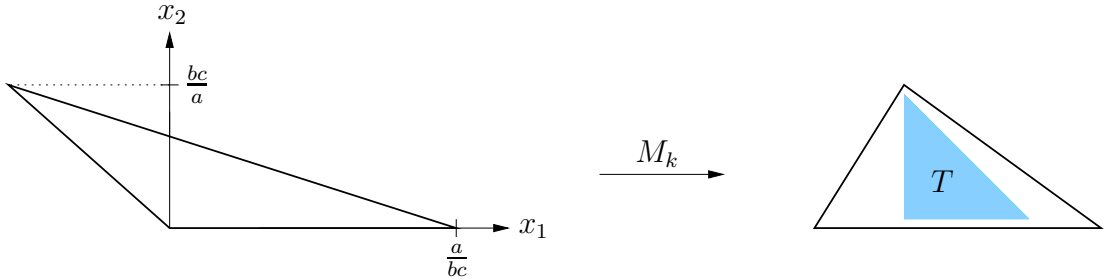}
 \end{center}
 \caption{The ball embedding providing the lower bound}
 \label{fig:shear}
\end{figure}

\m
By Proposition~\ref{prop:affwidth} and its proof, 
we can assume that the longest edge of~$\Delta_{\fm}$ is on the $x_1$-axis, starting in~$(0,0)$,
and we know that this longest length is~$\frac{a}{bc}$ and that the height over this edge is~$\frac{bc}a$,
see the left drawing in Figure~\ref{fig:shear}.
The claim for $(a,b,c) = (1,1,1)$ is clear. 

\begin{lemma} \label{le:alg}
Assume that $a \geqslant b \geqslant c$ and $(a,b,c) \neq (1,1,1)$. Then $\frac{a}{bc} > 2 \,\frac{bc}{a}$.
\end{lemma}

\proof 
The claim is true for $(a,b,c) = (2,1,1)$, and going down the Markov tree, 
after each mutation the longest edge becomes longer.
\proofend

By Lemma~\ref{le:alg}, after applying a matrix
$M_k = \begin{brsm} 1&k\\0&1\end{brsm}$ with $k \in \ZZ$, 
we can assume that the top vertex of~$\Delta_\fm$
lies over the interior of the horizontal side.
Fix a small $\varepsilon >0$.
Again by Lemma~\ref{le:alg},
the triangle $\Delta_\fm$ contains in its interior a right-angled triangle~$T$ 
of equal width and height $\frac{bc}a - \frac{\varepsilon}2$.
(If the top vertex of~$\Delta_\fm$ lies over the left half of the lower edge, 
we take~$T$ as in Figure~\ref{fig:shear}, if not, we take its reflection about the $x_2$-axis.) 
There exists a symplectic embedding 
$$
\varphi \colon B^4(\tfrac{bc}a - \varepsilon) \,\to\, T \times ]0,1[^2 \;\subset\; \RR^2(x) \times \RR^2(y) ,
$$
see \cite[Section~5]{Tra95} and \cite[Section~2]{Sch03} for two different constructions.
Therefore, $\varphi$ symplectically embeds 
$B^4(\tfrac{bc}a - \varepsilon)$ into 
$\Int \Delta_{\fm} \times \RR^2 / \ZZ^2 \subset \CP^2$. 
By Lemma~\ref{lem:slides} we can assume that $L_{a}^\fm \sqcup L_{b}^\fm \sqcup L_{c}^{\fm}$
is disjoint from the image of this embedding.
Since $\varepsilon >0$ was arbitrary, Step~1 is proved.

%

\b \noindent
{\bf Step 2.}
$w_{\G} (\CP^2 \setminus L_a^{\fm}) \leqslant \frac{bc}a$

\m 
Let $\gf \colon B^4(\ga) \to \CP^2 \setminus L_a^\fm$ be a symplectic embedding. 
We fix $\ga' < \ga$, restrict $\gf$ to the closed ball~$\overline B^4(\ga')$,
and set $B = \gf \bigl(\overline B^4(\ga') \bigr)$.
We show that $\ga' \leqslant \frac{bc}{a}$.

\begin{figure}[h] 
 \begin{center}
  \psfrag{B}{$\pi_\fm (B)$} \psfrag{H}{$H$} \psfrag{H'}{$H'$} \psfrag{V}{$V$}
  \leavevmode\epsfbox{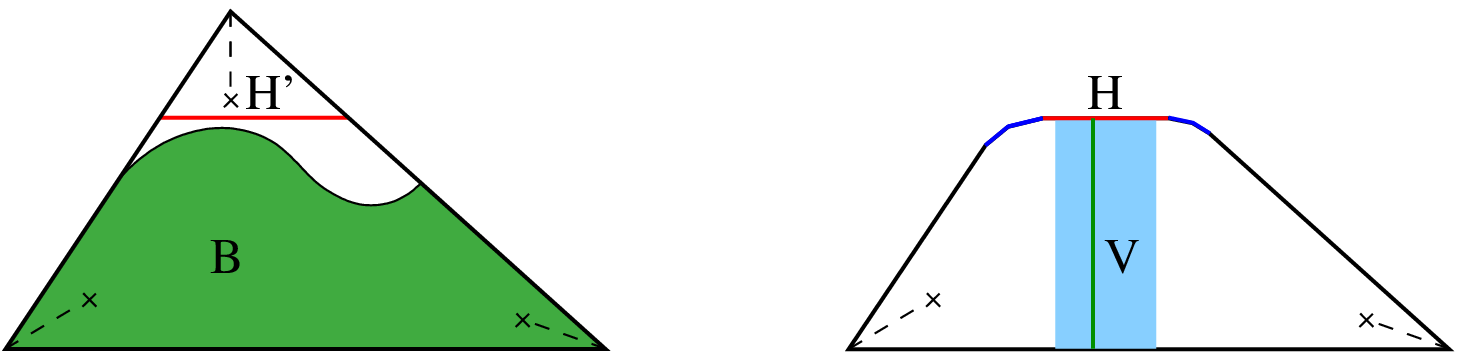}
 \end{center}
 \caption{From $\CP^2$ to $(M,\omega)$}
 \label{fig:blow-up}
\end{figure}

We work in the normal form of $\Delta_\fm$ used in the proof of Step~1 (cf.\ the right triangle in Figure~\ref{fig:shear}). 
The fibre of~$\pi_\fm$ over the top vertex is the core circle of the pinwheel~$L_a^\fm$.
Since $B$ and $L_a^\fm$ are disjoint and since $B$ is compact, there exists a neighbourhood of the top vertex in~$\Delta_\fm$
whose preimage is disjoint from~$B$.
Let $H'$ be a horizontal segment in~$\Delta_\fm$ strictly above~$\pi_\fm (B)$
as in the left drawing of Figure~\ref{fig:blow-up}.
By Lemma~\ref{lem:slides} we can assume that $\pi_{\fm}(L_a^\fm)$ lies strictly above~$H'$.
We now perform on~$\CP^2$ a symplectic cut along~$H'$, 
see~\cite{Le95}. 
Except for $a=2$, the resulting manifold is singular over the endpoints of~$H'$.
We resolve each singularity 
by applying finitely many symplectic cuts,  
so small that the ball~$B$ is untouched, and 
obtain a smooth closed symplectic manifold~$(M,\omega)$ fibering over a polygon 
as drawn on the right of Figure~\ref{fig:blow-up},
see for instance \cite[Section 4.5]{Eva21} for details.
The segment~$H'$ has now become a shorter segment $H$, over which lies a smooth symplectic
sphere~$S_H$.

%

The symplectic form $\omega$ over the rectangle on the RHS of Figure~\ref{fig:blow-up} 
is exactly the standard product form
on a sphere and an open cylinder. 
Take a vertical segment~$V$ from $H$ to the $x_1$-axis. 
Then there is an $\omega$-symplectic embedded sphere~$S_V$ over the vertical segment~$V$.
For $A := [S_V] \in H_2(M;\ZZ)$ we have $A \cdot A = 0$ and $c_1(M,\omega) (A) = 2$,
and $\langle [\omega], A \rangle = \int_{S_V}\omega = \frac{bc}{a} - \varepsilon =: h$.

\begin{proposition} \label{p:wG}
$w_{\G} (M,\omega) \leqslant h$.
\end{proposition}

Recall that by construction, $B = \varphi \bigl( \overline B^4(\alpha') \bigr)$ is a symplectic ball
in~$(M,\omega)$, whence $\ga' \leqslant w_{\G} (M,\omega)$.
The proposition thus implies that $\alpha' \leqslant h < \frac{bc}{a}$,
and therefore $w_{\G}(\CP^2 \setminus L_a^\fm) \leqslant \frac{bc}{a}$.

\medskip \noindent
{\it Proof of Proposition~\ref{p:wG}.}
We follow Gromov’s proof of his non-squeezing theorem and use Taubes’ Gromov invariant. 
Let $\cj (\omega)$ be the space of $\omega$-compatible almost complex structures 
on~$M$. 
For $J \in \cj (\omega)$ and a $J$-holomorphic curve~$C = u (\Sigma)$ in class~$A$,
the virtual genus is
$$
1 + \frac 12 \bigl( C \cdot C - c_1(C) \bigr) \,=\, 1 + \frac 12 \bigl( A^2 - c_1(A) \bigr) \,=\, 0
$$
and thus is equal to the genus. 
Since $A \cdot [S_H] = S_V \cdot S_H = 1$, the curve~$u$ is not a multiple covering. 
By the equality case in the adjunction formula we see that $C$ is embedded.
Furthermore, $c_1(A) = 2 \geqslant 1$. We conclude with the help of~\cite{HLS97} that $C$ is regular. 

Take a $J \in \cj (\omega)$ such that the symplectic sphere~$S_V$ over~$V$
is $J$-holomorphic.
(For instance, one can take the usual~$J$ on the product of a sphere and a cylinder 
over the rectangle on the RHS of Figure~\ref{fig:blow-up}, 
and extend this $J$ to an element of~$\cj (\omega)$.)
Fix a point $z_0$ over~$V$.  
Assume that $C$ is another closed $J$-holomorphic curve through~$z_0$
in class~$A$.
If $C \neq S_V$, then $C \cdot S_V \geqslant 1$ by positivity of intersection,
while $C \cdot S_V = A^2 =0$.
Hence $C = S_V$, and so Taubes' Gromov invariant $\Gr (A)$ is~$1$,
see \cite{Ta96} or~\cite[Section~1]{Mc97}.

Now let $\phi \colon B^4(\beta) \to (M,\omega)$ be a symplectic embedding.
After composing $\phi$ with a Hamiltonian diffeomorphism of $(M,\omega)$, 
we can assume that $\phi (0) = z_0$.
Take $J' \in \cj (\omega)$ such that $J' |_{\phi (B^4(\beta))} = \phi_* J_0$,
where $J_0$ is the standard complex structure~$i \oplus i$ on $B^4(\beta) \subset \CC^2$.
As we have seen, every $J'$-holomorphic sphere in class~$A$ is regular, and 
$\Gr (A) =1$.
Therefore, there exists a closed, possibly disconnected, $J'$-holomorphic curve 
$u \colon \Sigma \to M$ such that $[u(\Sigma)] = A$ and $u(\Sigma) \ni \phi (0)$.
Let $\Sigma'$ be a component of~$\Sigma$ with $u(\Sigma') \ni \phi (0)$.
Gromov's monotonicity argument from~\cite{Gr85} (see also~\cite{Hu97})
shows that $\beta \leqslant \int_{\Sigma'} u^* \omega$.
Since 
$$
\int_{\Sigma'} u^* \omega \,\leq\, \int_{\Sigma} u^* \omega \,=\, \langle [\omega], A \rangle
\,=\, h,
$$
we conclude that $\beta \leqslant h$.
\proofend

\begin{remarks}
{\rm
{\bf 1.}
One can give a different proof of Proposition~\ref{p:wG} by using the fact that the Seiberg--Witten invariant 
of~$A$ is~$1$ and Hutching's 
recent construction of a symplectic capacity in~\cite{Hu22} and his Lemma~7.4 therein.

\m \noindent
{\bf 2.}
The weaker statement that the Gromov width of $\CP^2 \setminus (L_{a}^\fm \sqcup L_{b}^\fm \sqcup L_{c}^{\fm})$
with $a \geqslant b,c$ is at most $\frac{bc}a$ follows more directly from the main result of~\cite{ChWo20}:
Given a symplectic embedding $B^4(\alpha) \to \CP^2 \setminus (L_{a}^\fm \sqcup L_{b}^\fm \sqcup L_{c}^{\fm})$
we apply the above symplectic cut procedure near each vertex of~$\Delta_\fm$
and obtain a symplectic embedding of a slightly smaller closed ball 
into a smooth closed toric symplectic manifold~$(M,\omega)$. 
For small enough cuts, 
the lattice width of the moment polytope of $(M,\omega)$ is arbitrarily close to the lattice width~$\frac{bc}{a}$
of~$\Delta_\fm$, 
whence $w_{\G}(M,\omega) \leqslant \frac{bc}{a}$ by \cite[Corollary~4.19]{ChWo20},
and hence also $w_{\G} \bigl(\CP^2 \setminus (L_{a}^\fm \sqcup L_{b}^\fm \sqcup L_{c}^{\fm}) \bigr) \leqslant \frac{bc}{a}$.

\m \noindent
{\bf 3.}
For a symplectic orbifold $(M,\omega)$ define the Gromov width as the supremum over~$\alpha$
for which $B^{2n}(\alpha)$ symplectically embeds into the smooth locus of~$M$.
The weighted projective plane $\CP (a^2,b^2,c^2)$ has (at most) three singularities $s_a, s_b, s_c$,
and the complement carries an effective Hamiltonian $T^2$-action,
whose moment polygon is~$\Delta_{(a^2,b^2,c^2)}$, for the suitably normalized symplectic form.
It is shown in~\cite{Bre21} that  
$\CP (a^2,b^2,c^2) \setminus \{s_a, s_b, s_c\}$ and $\CP^2 \setminus (L_{a}^\fm \sqcup L_{b}^\fm \sqcup L_{c}^{\fm})$
are symplectomorphic.
The previous remark and the construction in Step~1 now show that $w_{\G} (\CP(a^2,b^2,c^2)) = \frac{bc}a$.

\m \noindent
{\bf 4.\ Relation to Biran's argument.}
In \cite{Bi01} Biran proved $w_{\G}(\CP^2 \setminus \RP^2) \leqslant \frac 12$
in the following way.
Take a smooth quadric $Q \subset \CP^2$, i.e.\ a smoothly embedded
symplectic sphere that is homologous to twice the class of a line.
The complement of the symplectic disc normal bundle~$DQ$ with discs of area~$\frac 12$ is~$\RP^2$.
(Starting from $Q = \{ z_0^2+z_1^2+z_2^2=0 \}$ Biran obtains as complement the standard~$\RP^2$, 
and in general one obtains a set Hamiltonian isotopic to this $\RP^2$.)
An adaptation of Gromov's non-squeezing theorem then yields
$$
w_{\G}(\CP^2 \setminus \RP^2) \,=\, w_{\G}(DQ) \,\leq\, \tfrac 12 .
$$

\begin{figure}[h] 
 \begin{center}
  \psfrag{Q}{$Q$} \psfrag{Q0}{$Q_0$} \psfrag{R}{$\RP^2$}
  \leavevmode\epsfbox{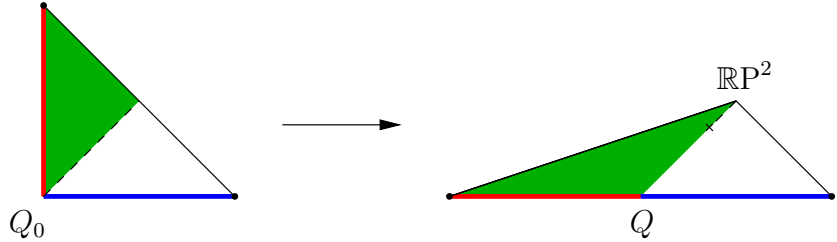}
 \end{center}
 \caption{The almost toric picture for Biran's decomposition $\CP^2 = DQ \sqcup \RP^2$}
 \label{fig:R3}
\end{figure}

%

The union of two lines $Q_0 = \{z_1 z_2 = 0\} \subset \CP^2$ is a singular quadric.
After a nodal trade, $Q_0$ becomes a smooth quadric~$Q \subset \CP^2$, 
``opposite" to~$\RP^2$, just as in Biran's decomposition.

\begin{question*}
Can the spaces $\CP^2 \setminus L_p^{\fm}$
for $p \geqslant 5$ be seen as symplectic disc bundles over (singular) symplectic surfaces?
\end{question*}
}
\end{remarks}

\appendix

\section{Pinwheels are ``monotone"}
\label{app:monotone}

Let $L_p \subset \CP^2$ be a Lagrangian pinwheel. 
The long exact sequences of homotopy and homology groups 
and the Hurewicz homomorphisms yield the commutative diagram
\begin{equation*} 
\xymatrixcolsep{3pc}
\xymatrix{ 
 \pi_2(L_p)  \ar@{->}[r] \ar@{->}[d]^{\alpha} & \pi_2(\CP^2)  \ar@{->}[r] \ar@{->}[d]^{\beta} &
 \pi_2(\CP^2, L_p)  \ar@{->}[r] \ar@{->}[d]^{\gamma} & \pi_1(L_p)  \ar@{->}[r] \ar@{->}[d]^{\delta} & 
0 \ar@{->}[d] 
\\
H_2(L_p)  \ar@{->}[r]  & H_2(\CP^2)  \ar@{->}[r] & H_2(\CP^2, L_p)  \ar@{->}[r] & H_1(L_p) \ar@{->}[r]  &
0 
}
\end{equation*}
Since $H_2(L_p) = 0$ and $\beta$ and $\delta$ are isomorphisms, 
the Five lemma implies that $\gamma$ is an isomorphism.
We thus have the short exact sequence
\begin{equation*} 
\xymatrixcolsep{3pc}
\xymatrix{ 
 0  \ar@{->}[r] & \pi_2(\CP^2) = \ZZ  \ar@{->}[r]^{\iota} &
 \pi_2(\CP^2, L_p)  \ar@{->}[r] & \ZZ_p = \pi_1(L_p)  \ar@{->}[r] & 
0 
}
\end{equation*}
In particular, $\pi_2(\CP^2,L_p) \cong \ZZ \oplus T$, 
where $T$ is a finite abelian group. 
(One can show that $T \cong \ZZ_{p'}$ for a divisor $p'$ of~$p$, see Lemma~2.16~(c) in~\cite{ES18}.)

Recall that a closed Lagrangian submanifold $L$ of a symplectic manifold~$(M,\omega)$ 
is called \emph{monotone}\/
if the area homomorphism~$\sigma_\omega$ and the Maslov homomorphism~$\sigma_\mu$
are positively proportional on~$\pi_2(M,L)$. 
For our immersed Lagrangian~$L_p$ the area homomorphism can be defined as for embedded Lagrangians:
$$
\sigma_\omega (A) \,:=\, \int_\DD u^* \omega_{\FS}
$$
where $u \colon (\DD, \partial \DD) \to (\CP^2, L_p)$ is any smooth map representing the class~$A$.

Assume we can meaningfully define a Maslov homomorphism $\sigma_\mu \colon \pi_2(\CP^2,L_p) \to \RR$
with the property that 
\begin{equation} \label{e:mumu}
\sigma_{\mu}(\iota ([\CP^1])) = \sigma_\mu ([\CP^1]) ,
\end{equation}
where $\sigma_\mu ([\CP^1]) = 2 c_1 ([\CP^1]) = 6$
denotes the usual Maslov index of the ``disc"~$\CP^1$ in~$\CP^2$.
Taking values in~$\RR$, the homomorphisms $\sigma_\omega$ and $\sigma_\mu$
on $\pi_2(\CP^2, L_p) \cong \ZZ \oplus T$ vanish on~$T$.
Let $(k,t) \in \ZZ \oplus T$ be the element corresponding to $\iota ([\CP^1])$
under this isomorphism. 
Then
$$
\begin{array}{lclclcl}
\sigma_\omega (k,0) &=& \sigma_{\omega}(k,t) &=& \sigma_{\omega}(\iota ([\CP^1])) &=& 1, \\ [0.4em]
\sigma_\mu (k,0) &=& \sigma_{\mu}(k,t) &=& \sigma_{\mu}(\iota ([\CP^1])) &=& 6.
\end{array}
$$
Hence $\sigma_\mu = 6 \, \sigma_{\omega}$.
We can thus rightfully say that $L_p \subset \CP^2$ is monotone.

It remains to define $\sigma_\mu$ on~$\pi_2(\CP^2,L_p)$ satisfying~\eqref{e:mumu}.
We present an element of $\pi_2(\CP^2, L_p)$ by a smooth map 
$u \colon (\DD, \partial \DD) \to (\CP^2, C)$,
where $C \subset L_p$ is the core circle.
Consider the map $u^p (z) = u(z^p)$, 
that represents $p [u]$.
Choose one of the $p$ branches of~$L_p$ at~$u^p(1)$.
Let $P_0$ be the Lagrangian plane at $u^p(1)$ tangent to this branch.
Following this branch along the path $v(t) := u^p(e^{2\pi i t})$, 
we obtain a path of Lagrangian planes $P_t \subset T_{v(t)} \CP^2$, $t \in [0,1]$.
Since the class of the loop~$v$ vanishes in~$\pi_1(C)$, the path~$P_t$ is a loop.
If we define $P_t$ by starting on another branch, we get the same loop 
(with a shift in the parametrization).
As in the definition for embedded Lagrangians, we now
choose a symplectic trivialization of $T\CP^2 |_{u^p(\DD)}$,
and define $\sigma_\mu ([u^p]) \in \ZZ$ as the Maslov index of~$P_t$
with respect to this trivialization.
Finally, we set $\sigma_\mu ([u]) = \frac 1p \2 \sigma_\mu ([u^p]) \in \QQ$.
Since $\CP^1$ can be represented by a map $u \colon (\DD, \partial \DD) \to (\CP^2, C)$
taking $\partial \DD$ to a point in~$C$,
this definition satisfies~\eqref{e:mumu}.

\end{document}